\documentclass[11pt]{article}%
\usepackage{amsmath}
\usepackage{amsfonts}
\usepackage{amssymb}
\usepackage{graphicx}
\usepackage{a4wide}
\newtheorem{theorem}{Theorem}
\newtheorem{proposition}[theorem]{Proposition}
\begin{document}

\title{{\Large On the separation theorems for convex sets on the unit sphere}}
\author{Constantin Z\u{a}linescu\thanks{Octav Mayer Institute of Mathematics, Ia\c{s}i
Branch of Romanian Academy, Ia\c{s}i, Romania, and University
\textquotedblleft Alexandru Ioan Cuza" Ia\c{s}i, Romania; email:
\texttt{zalinesc@uaic.ro}.}}
\date{}
\maketitle

\begin{abstract}
In this short note we provide a new proof of the recent result of Han and
Nashimura on the separation of spherical convex sets established in
arXiv:2002.06558. Our proof is based on a result stated in locally convex spaces.

\end{abstract}

Recall that the support function of the nonempty subset $A$ of the locally
convex space $X$ is the function
\[
\sigma_{A}:X^{\ast}\rightarrow\overline{\mathbb{R}},\quad\sigma_{A}(x^{\ast
}):=\sup\{\left\langle x,x^{\ast}\right\rangle \mid x\in A\}\quad(x^{\ast}\in
X^{\ast}),
\]
where $X^{\ast}$ is the topological dual of $X$ and $\left\langle x,x^{\ast
}\right\rangle :=x^{\ast}(x)$ for $x\in X$ and $x^{\ast}\in X^{\ast}$. It is
clear that $\sigma_{A}$ is convex, positively homogeneous (hence sublinear),
$w^{\ast}$-lsc and $\sigma_{A}=\sigma_{\overline{\operatorname*{conv}}A}$; in
particular, $[\sigma_{A}\leq\alpha]$ $(:=\{x^{\ast}\in X^{\ast}\mid\sigma
_{A}(x^{\ast})\leq\alpha\})$ is $w^{\ast}$-closed\ for every $\alpha
\in\mathbb{R}$. Moreover, $\partial\sigma_{A}(0)=\overline
{\operatorname*{conv}}A$, where $\operatorname*{conv}A$ (resp.\ $\overline
{\operatorname*{conv}}A$) denotes the convex (resp.\ the closed convex) hull
of $A$ and
\[
\partial\sigma_{A}(0):=\{x\in X\mid\left\langle x,x^{\ast}\right\rangle
\leq\sigma_{A}(x^{\ast})\quad\forall x^{\ast}\in X^{\ast}\}.
\]

\begin{proposition}
\label{p-hn0}Let $X$ be a Hausdorff locally convex space, $A\subseteq X$ a
nonempty set and $\alpha\in\mathbb{R}$. Set $D_{\alpha}:=\{x^{\ast}\in
X^{\ast}\mid\left\langle x,x^{\ast}\right\rangle <\alpha\ \forall x\in A\}.$

\emph{(i)} If $A$ is open, then $[\sigma_{A}\leq\alpha]\setminus\{o\}\subseteq
D_{\alpha}\subseteq\lbrack\sigma_{A}\leq\alpha],$ where $[\sigma_{A}\leq
\alpha]:=\{x^{\ast}\in X^{\ast}\mid\sigma_{A}(x^{\ast})\leq\alpha\}$ and $o$
is the zero element of $X^{\ast}$; consequently, $D_{\alpha}\cup\{o\}$ is
$w^{\ast}$-closed.

\emph{(ii)} If $A$ is $w$-compact, then $D_{\alpha}$ is
$\sigma$-open, where $\sigma$ is the Mackey topology on $X^{\ast}$
(cf.\ \cite[p.\ 260]{Koe69}).
\end{proposition}

Proof. (i) Consider $x^{\ast}\in\lbrack\sigma_{A}\leq\alpha]\setminus\{o\}$.
Assume that $x^{\ast}\notin D_{\alpha}$; then there exists $x_{0}\in A$ such
that $\alpha\leq\left\langle x_{0},x^{\ast}\right\rangle $ $(\leq\alpha),$ and
so $\left\langle x_{0},x^{\ast}\right\rangle =\alpha$. It follows that
$\left\langle x-x_{0},x^{\ast}\right\rangle =\left\langle x,x^{\ast
}\right\rangle -\alpha\leq0$ for every $x\in A$, and so $\left\langle
x^{\prime},x^{\ast}\right\rangle \leq0$ for every $x^{\prime}\in A-x_{0}=:V$.
Since $A$ is open and $x_{0}\in A$, $V$ is a neighborhood of the origin of
$X$, whence $V$ is absorbing; consequently $\left\langle x^{\prime},x^{\ast
}\right\rangle =0$ for every $x^{\prime}\in X$, whence the contradiction
$x^{\ast}=o$. Hence, $[\sigma_{A}\leq\alpha]\setminus\{o\}\subseteq D_{\alpha
}$ holds. The second inclusion is obvious.

Consequently, if $0=\sigma_{A}(o)\leq\alpha,$ then $D_{\alpha}\cup\left\{
o\right\}  =[\sigma_{A}\leq\alpha]$ is $w^{\ast}$-closed because $\sigma_{A}$
is $w^{\ast}$-lsc. If $\alpha<0,$ then $D_{\alpha}=[\sigma_{A}\leq\alpha],$
and so $D_{\alpha}\cup\left\{  o\right\}  $ is $w^{\ast}$-closed because
$w^{\ast}$ is Hausdorff.

(ii) Take $x_{0}^{\ast}\in D_{\alpha}$; then there exists $x_{0}\in A$ such
that $\left\langle x,x_{0}^{\ast}\right\rangle \leq\left\langle x_{0}%
,x_{0}^{\ast}\right\rangle $ $(<\alpha)$ for every $x\in A$. Then
$V:=\{x^{\ast}\in X^{\ast}\mid\left\vert \left\langle x,x^{\ast}\right\rangle
\right\vert \leq\gamma\ \forall x\in A\}=\gamma([-1,1]A)^{\circ}$ is a
$\sigma$-neighborhood of $o$, where $\gamma:=\tfrac{1}{2}(\alpha-\left\langle
x_{0},x_{0}^{\ast}\right\rangle )>0$. Then, for $x^{\ast}\in V$ and $x\in A$
one has
\[
\left\langle x,x_{0}^{\ast}+x^{\ast}\right\rangle =\left\langle x,x_{0}^{\ast
}\right\rangle +\left\langle x,x^{\ast}\right\rangle \leq\left\langle
x_{0},x_{0}^{\ast}\right\rangle +\left\langle x,x^{\ast}\right\rangle
\leq\left\langle x_{0},x_{0}^{\ast}\right\rangle +\gamma=\tfrac{1}{2}%
(\alpha+\left\langle x_{0},x_{0}^{\ast}\right\rangle )<\alpha;
\]
consequently, $x_{0}^{\ast}+V\subseteq D_{\alpha}$. Therefore,
$D_{\alpha}$ is $\sigma$-open. \hfill $\square$

\medskip

In the next result we apply Proposition \ref{p-hn0} for
$X:=\mathbb{R}^{n}$ $(n\geq2)$ endowed with the usual inner product
$\left\langle \cdot ,\cdot\right\rangle $ and its Euclidean norm
$\left\Vert \cdot\right\Vert $, $X^{\ast}$ being identified with
$\mathbb{R}^{n}$; the unit sphere is denoted, as usual, by
$\mathbb{S}^{n-1}$.

\begin{theorem}
\label{t-HN} Let $B_{1},B_{2}\subseteq\mathbb{S}^{n-1}$ $(n\geq2)$ be nonempty
sets such that $P_{1}:=\mathbb{P}B_{1}$ and $P_{2}:=\mathbb{P}B_{2}$ are
convex (equivalently, $B_{1}$ and $B_{2}$ are spherically convex), where
$\mathbb{P}:={}]0,\infty\lbrack$. Consider the set $E:=E_{1}\cap E_{2},$
where
\begin{equation}
E_{1}:=\{u\in\mathbb{S}^{n-1}\mid\left\langle b_{1},u\right\rangle >0\ \forall
b_{1}\in B_{1}\},\ \ E_{2}:=\{u\in\mathbb{S}^{n-1}\mid\left\langle
b_{2},u\right\rangle <0\ \forall b_{2}\in B_{2}\}. \label{r-hno}%
\end{equation}
Then the following assertions hold:

\emph{(i)} If $E\neq\emptyset,$ then $B_{1}\cap B_{2}=\emptyset$.

\emph{(ii)} If $B_{1}\cap B_{2}=\emptyset$ and $P_{1}$, $P_{2}$ are open, then
$E$ is nonempty and $\mathbb{R}_{+}E$ is a pointed, closed and convex cone;
consequently, $E$ is a closed spherical convex set.

\emph{(iii)} If $B_{1}\cap B_{2}=\emptyset$ and $B_{1}$, $B_{2}$ are closed,
then $E$ is nonempty and $\mathbb{P}E$ is convex and open; consequently, $E$
is an open spherical convex set.
\end{theorem}

Proof. Clearly,
\begin{gather}
\mathbb{P}E_{1}=\{u\in\mathbb{R}^{n}\mid\left\langle x_{1},u\right\rangle
>0\ \forall x_{1}\in P_{1}\}=\{u\in\mathbb{R}^{n}\mid\left\langle
x_{1},u\right\rangle >0\ \forall x_{1}\in B_{1}\},\label{r-hn11}\\
\mathbb{P}E_{2}=\{u\in\mathbb{R}^{n}\mid\left\langle x_{2},u\right\rangle
<0\ \forall x_{2}\in P_{2}\}=\{u\in\mathbb{R}^{n}\mid\left\langle
x_{2},u\right\rangle <0\ \forall x_{2}\in B_{2}\}, \label{r-hn12}%
\end{gather}
and so $\mathbb{P}E_{1}$, $\mathbb{P}E_{2}$ are convex sets. Because
$B_{1},B_{2},E_{1},E_{2}\subset\mathbb{S}^{n-1},$ one has $P_{1}\cap
P_{2}=\mathbb{P}(B_{1}\cap B_{2})$ and $\mathbb{P}E=\mathbb{P}E_{1}%
\cap\mathbb{P}E_{2}$; hence $\mathbb{P}E$ is also convex.

(i) The assertion is obviously true.

(ii) Because $B_{1}\cap B_{2}=\emptyset$, one has $P_{1}\cap P_{2}=\emptyset$.
As $P_{1}$ and $P_{2}$ are nonempty disjoint open convex sets, by a well-known
separation theorem, there exists $u\in\mathbb{S}^{n-1}$ and $\gamma
\in\mathbb{R}$ such that $\left\langle x_{1},u\right\rangle \geq\gamma
\geq\left\langle x_{2},u\right\rangle $ for all $(x_{1},x_{2})\in P_{1}\times
P_{2}$. It follows that $\gamma=0$ because $P_{1}=\mathbb{P}P_{1}$ and
$P_{2}=\mathbb{P}P_{2};$ moreover, because $P_{1}$ and $P_{2}$ are open, one
obtains that $\left\langle x_{1},u\right\rangle >0>\left\langle x_{2}%
,u\right\rangle $ for $(x_{1},x_{2})\in P_{1}\times P_{2}$. As $B_{1}\subseteq
P_{1}$ and $B_{2}\subseteq P_{2},$ one has $u\in E_{1}\cap E_{2}$, and so
$E\neq\emptyset$. Applying Proposition \ref{p-hn0}(i) for $P_{2}$ and
$\alpha=0$, and taking into account the first expression of $\mathbb{P}E_{2}$
in (\ref{r-hn12}), one obtains that $\mathbb{R}_{+}E_{2}$ $(=(\mathbb{P}%
E_{2})\cup\{o\})$ is closed; the closedness of $\mathbb{R}_{+}E_{1}$ is
obtained similarly, and so $\mathbb{R}_{+}E$ is closed.

(iii) Because for $i\in\{1,2\}$ the set $B_{i}$ $(\subset\mathbb{S}^{n-1})$ is
closed, $B_{i}$ is compact. Using the second expression of $\mathbb{P}E_{2}$
in (\ref{r-hn12}) and Proposition \ref{p-hn0}(ii), the set $\mathbb{P}E_{2}$
is open; the openness of $\mathbb{P}E_{1}$ is obtained similarly, and so
$\mathbb{P}E$ is open.

It remains to show that $E$ is nonempty. For this aim, observe first that
$C_{i}:=\operatorname*{conv}B_{i}$ $(\subseteq P_{i})$ is compact and
$P_{i}=\mathbb{P}C_{i}$ for $i\in\{1,2\};$ in particular, $C_{1}\cap
C_{2}=\emptyset$ and $o\notin C_{1}\cup C_{2}$. We claim that there exists
$r\in\mathbb{P}$ such that $\left[  \mathbb{P}(B_{1}+rU)\right]  \cap\left[
\mathbb{P}(B_{2}+rU)\right]  =\emptyset$, where $U$ is the unit closed ball of
$\mathbb{R}^{n}$. In the contrary case, having the sequence $(r_{n})_{n\geq
1}\subset\mathbb{P}$ with $r_{n}\rightarrow0$, for every $n\in\mathbb{N}$ and
every $i\in\{1,2\}$ there exist $t_{n}^{i}\in\mathbb{P}$, $v_{n}^{i}\in C_{i}$
and $x_{n}^{i}\in U$ such that $t_{n}^{1}(v_{n}^{1}+r_{n}x_{n}^{1})=t_{n}%
^{2}(v_{n}^{2}+r_{n}x_{n}^{2})$; clearly, we may take $t_{n}^{1}=1$ for
$n\in\mathbb{N}$. Passing to subsequences, we may assume that $t_{n}%
^{2}\rightarrow t\in\lbrack0,\infty],$ $v_{n}^{i}\rightarrow v_{i}%
\in\operatorname*{cl}C_{i}=C_{i}$ for $i\in\{1,2\}$. If $t=0$ we get the
contradiction $C_{1}\ni v_{1}=o;$ if $t=\infty$ we get the contradiction
$C_{2}\ni v_{2}=o$. Hence $t\in\mathbb{P}$, whence $v_{1}=tv_{2}\in
\mathbb{P}C_{1}\cap\mathbb{P}C_{2}=\emptyset.\ $Therefore, our claim is true,
that is, there exists $r\in\mathbb{P}$ such that $\left[  \mathbb{P}%
(C_{1}+rU)\right]  \cap\left[  \mathbb{P}(C_{2}+rU)\right]  =\emptyset$.
Because the sets $\mathbb{P}(C_{i}+rU)$ with $i\in\{1,2\}$ are convex sets
(even with nonempty interior), there exist $u\in\mathbb{S}^{n-1}$ and
$\alpha\in\mathbb{R}$ such that
\[
\left\langle t_{1}(v_{1}+rx_{1}),u\right\rangle \geq\alpha\geq\left\langle
t_{2}(v_{2}+rx_{2}),u\right\rangle \quad\forall t_{1},t_{2}\in\mathbb{P}%
,\ x_{1},x_{2}\in U,\ v_{1}\in C_{1},\ v_{2}\in C_{2}.
\]
It follows that $\alpha=0$ and $\left\langle v_{1},u\right\rangle
-r\geq 0\geq\left\langle v_{2},u\right\rangle +r,$ whence $u\in
E_{1}\cap E_{2}=E$.  \hfill $\square$

\medskip Theorem \ref{t-HN} is established by Han and Nashimura in
\cite[Th.\ 1]{HanNis:20}. Note that from Theorem \ref{t-HN}(iii)
one obtains the strict separation theorem of disjoint closed
strongly convex sets on a sphere by a hyperplane through the centre
from \cite[Th.\ 1]{Bak:67} (attributed to an unpublished paper by BC
Rennie), while from Theorem \ref{t-HN}(ii) one obtains the
separation theorem of convex bodies with disjoint interiors stated
in \cite[Lem.\ 2]{Las:15}.

\end{document}